\documentclass[12pt]{article}
\usepackage{MnSymbol}
\usepackage{amsthm}

\DeclareMathAlphabet\mathbb{U}{msb}{m}{n}
\DeclareMathAlphabet\Bbb{U}{msb}{m}{n}

\usepackage{vmargin}
\usepackage{wasysym}
\usepackage[active]{srcltx}
\usepackage[stable]{footmisc}
\usepackage{caption}
\usepackage{pdfpages}
\usepackage[matrix, arrow,all,cmtip,color]{xy}
\usepackage{url}
\usepackage{epigraph}

\newcommand{\bi}{\begin{itemize}}
\newcommand{\ei}{\end{itemize}}
\newcommand{\bd}{\begin{description}}
\newcommand{\ed}{\end{description}}

\theoremstyle{plain}

\theoremstyle{definition}
\newtheorem{defn}{Definition}

\def\bqqq{\begin{quote}}
\def\eqqq{\end{quote}}

\newcommand*{\longhookrightarrow}{\ensuremath{\lhook\joinrel\relbar\joinrel\rightarrow}}

\def\lra{\longrightarrow}

\def\rtt{\,\rightthreetimes\,}

\def\rrt#1#2#3#4#5#6{\xymatrix{ {#1} \ar[r]^{} \ar@{->}[d]_{#2} & {#4} \ar[d]^{#5} \\ {#3}  \ar[r] \ar@{-->}[ur]^{}& {#6} }}

\def\union{\cup}

\def\RR{\Bbb R}

\begin{document}
\title{
%%%%n 
Separation axioms as lifting properties
}
\author{
%%%%n 
misha gavrilovich\thanks{A draft; comments welcome. {\tt mi\!\!\!ishap\!\!\!p@sd\!\!\!df.org http://mishap.sdf.org.} }
}
\maketitle

\begin{abstract}
We observe that many of the separation axioms of topology (including $T_0-T_4$)
can be expressed concisely and uniformly
in terms of category theory as lifting properties (in the sense of 
Quillen model categories) with respect to (usually open) continuous maps
of finite spaces (involving up to \ensuremath{4} points)
and the real line. 
\end{abstract}

\section{Introduction}

We observe that separation axioms of topology including $T_0-T_4$ 
can be expressed concisely and uniformly 
in terms of category theory as Quillen lifting properties with respect to (usually open) continuous maps 
of finite spaces (involving up to \ensuremath{4} points) 
and the real line. To make the exposition as self-contained as possible, 
we took the Wikipedia page on the separation axioms and added there 
the lifting properties formulae expressing what is said there in words.

No attempt is made here to explore the expressive power of the Quillen lifting property; 
for example, 
the note leaves out a reformulation of compactness
in terms of the Quillen lifting property and maps of finite topological spaces, as well as iterated lifting properties. 
See [G] for an attempt to suggest a context for these observations;
in particular, [G] discusses this and other examples, e.g. 
a finite group being nilpotent, solvable, p-group, 
a module being projective, injective, and a map being injective, surjective. 

{\bf Acknowledgements.} Thanks are due to the authors of the Wikipedia page 
on the lifting properties. I thank S.V.Ivanov for suggestions which helped to improve the exposition. See [G, DMG] for more.  

\section{Preliminaries}

\subsection{Quillen lifting property}

\begin{defn} Let $f$ and $g$ be a pair of morphisms from a category. 
We say that `` $f$ has the left lifting property wrt $g$ ",
`` $f$ is (left) orthogonal to $g$ ",
and write 
$f \,\rightthreetimes\, g$, iff for each $i:A\longrightarrow X$, $j:B\longrightarrow Y$ such that $ig=fj$ (``the square commutes"),
there is $j':B\longrightarrow X$ such that $fj'=i$ and $j'g=j$ (``there is a diagonal
making the diagram commute").
\end{defn}
\def\rrt#1#2#3#4#5#6{\xymatrix{ {#1} \ar[r]^{} \ar@{->}[d]_{#2} & {#4} \ar[d]^{#5} \\ {#3}  \ar[r] \ar@{-->}[ur]^{}& {#6} }}
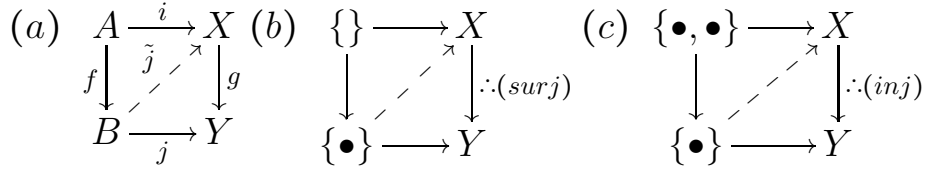
\begin{figure}
\begin{center}
\large
$ (a)\ \xymatrix{ A \ar[r]^{i} \ar@{->}[d]_f & X \ar[d]^g \\ B \ar[r]_-{j} \ar@{-->}[ur]^{{\tilde j}}& Y }$% \ 
%$\rrt ABXY$\ \ \ \
$(b)\  \rrt  {\{\}}  {} {\{\bullet\}}  X {\therefore(surj)} Y $%\  
$(c)\  \rrt {\{\bullet,\bullet\}} {} {\{\bullet\}}  X {\therefore(inj)} Y $%\ \ 
%$(d)\  \rrt X {\therefore(inj)} {Y} {\{x,y\}} {} {\{x=y\}}$\ 
\end{center}
\caption{\label{fig1}\normalsize
Lifting properties. %Dots $\therefore$ indicate free variables and what property of these variables is being defined;
%, i.e.~a property of what is being defined and how is it to be labelled 
%in a diagram chasing calculation, ``$\therefore(surj)$" reads as: 
%given a (valid) diagram, add label $(surj)$ to the corresponding arrow.\newline
 (a) The definition of a lifting property $f\rtt g$. 
%: for each $i:A\lra X$ and $j:B\lra Y$
%making the square commutative, i.e.~$f\circ j=i\circ g$, there is a diagonal arrow $\tilde j:B\lra X$ making the total diagram
%$A\xra f B\xra {\tilde j} X\xra g Y, A\xra i X, B\xra j Y$ commutative, i.e.~$f\circ \tilde j=i$ and $\tilde j\circ g=j$.
 (b) $X\lra Y$ is surjective %\newline
 (c) $X\lra Y$ is injective 
%; $X\lra Y$ is an epicmorphism if we forget %never use
%that $\{\bullet\}$ denotes a singleton (rather than an arbitrary object
%and thus $\{\bullet,\bullet\}\lra\{\bullet\}$ denotes an arbitrary morphism $Z\sqcup Z\xra{(id,id)} Z$)\newline
% (d) $X\lra Y$ is injective, in the category of Sets; $\pi_0(X)\lra\pi_0(Y)$ is injective if 
%  the diagram is interpreted in the category
%of topological spaces.
}
\end{figure}

A useful intuition is to think that the property of left-lifting against each map in a
class \ensuremath{C} is a kind of negation of the property of being in \ensuremath{C}, and that
right-lifting is another kind of negation. 
For example, the Sierpinski space $S$  consisting of one open point and one closed point, 
is perhaps the simplest counterexample
to the separation axiom $T_1$, and a space $X$ satisfies $T_1$ iff  
$$ S\longrightarrow {pt} \,\rightthreetimes\,  \ensuremath{X} \longrightarrow {pt}$$
where $S\longrightarrow {pt}$, resp.~$X\longrightarrow {pt}$, denote the map sending $S$, resp.~$X$, into the space
consisting of a single point.

\subsection{Notation for maps of finite topological spaces}

A topological space comes with a {\em specialisation preorder} on its points: for
points $x,y \in X$,  $x \leq y$ iff $y \in cl x$ , or equivalently,  a category whose
objects are points of \ensuremath{X} and there is a unique morphism $x{\small\searrow}y$ iff $y \in cl x$.

For a finite topological space X, the specialisation preorder or equivalently
the category uniquely determines the space: a subset of \ensuremath{X} is closed iff it is
downward closed, or equivalently, there are no morphisms going outside the
subset.

The monotone maps (i.e. functors) are the continuous maps for this topology.

We denote a finite topological space by a list of the arrows (morphisms) in
the corresponding category; '$\leftrightarrow $' denotes an isomorphism and '$=$' denotes
the identity morphism.  An arrow between two such lists denotes
a continuous map (a functor) which sends each point to the correspondingly
labelled point, but possibly turning some morphisms into identity morphisms,
thus gluing some points. 

Thus, each point goes to "itself" and  $$
     \{a,b\}\longrightarrow \{a{\small\searrow}b\}\longrightarrow \{a\leftrightarrow b\}\longrightarrow \{a=b\}
$$
denotes
$$
   (discrete\ space\ on\ two\ points)\longrightarrow (Sierpinski\ space)\longrightarrow (antidiscrete\ space)\longrightarrow (single\ point)
$$

In $A \longrightarrow  B$, each object and each morphism in $A$ necessarily appears in \ensuremath{B} as well; we avoid listing 
the same object or morphism twice. Thus 
both 
$$
\{a\} \longrightarrow  \{a,b\}\text{ and } \{a\} \longrightarrow  \{b\}
$$  denote the same map from a single point to the discrete space with two points.
Both 
 $$\{a{\small\swarrow}U{\small\searrow}x{\small\swarrow}V{\small\searrow}b\}\longrightarrow \{a{\small\swarrow}U=x=V{\small\searrow}b\}\text{ and }\{a{\small\swarrow}U{\small\searrow}x{\small\swarrow}V{\small\searrow}b\}\longrightarrow \{U=x=V\}$$
denote the morphism gluing points $U,x,V$.

In $\{a{\small\searrow}b\}$, the point $a$ is open and point \ensuremath{b} is closed.

\section{Separation Axioms}

Let \ensuremath{X} be a topological space. Then two points \ensuremath{x} and \ensuremath{y} in \ensuremath{X} are {\em topologically distinguishable}
iff the map $\{x\leftrightarrow y\} \longrightarrow  X$ is not continuous, i.e. 
iff %if they do not have exactly the same neighbourhoods (or equivalently the same open neighbourhoods); 
at least one of them has an open neighbourhood which is not a neighbourhood of the other.

Two points \ensuremath{x} and \ensuremath{y} are {\em separated} iff neither $ \{x{\small\searrow}y\} \longrightarrow  X$ nor $\{x{\small\searrow}y\} \longrightarrow  X$ is continuous, 
i.e~each of them has a neighbourhood that is not a neighbourhood of the other; 
in other words, neither belongs to the other's closure, $x \notin cl\, x$ and $y \notin cl\, x$. 
More generally, two subsets A and \ensuremath{B} of \ensuremath{X} are {\em separated} iff each is disjoint from the other's closure, 
i.e.~$ A\cap cl \ensuremath{B} = B\cap cl A = \emptyset $. 
(The closures themselves do not have to be disjoint.) In other words, the map
$ i_{AB} : \ensuremath{X} \longrightarrow  \{A \leftrightarrow  \ensuremath{x} \leftrightarrow  B\}$
sending the subset $A$ to the point $A$, the subset $B$ to the point $B$, and the rest to the point $x$,  
factors both as 
 $$ 
X \longrightarrow   \{A \leftrightarrow  U_A {\small\searrow} \ensuremath{x} \leftrightarrow  B\} \longrightarrow  \{A \leftrightarrow  U_A = \ensuremath{x} \leftrightarrow  B\} $$
and 
 $$ 
X \longrightarrow   \{A \leftrightarrow  \ensuremath{x} {\small\swarrow} U_B \leftrightarrow  B\} \longrightarrow  \{A \leftrightarrow  \ensuremath{x} = U_B \leftrightarrow  B\} $$ 
here the preimage of $x,B$, resp. $x,A$ is a closed subset containing $B$, resp. $A$, and disjoint from $A$, resp. $B$.
All of the remaining conditions for separation of sets may also be applied to points (or to a point and a set) 
by using singleton sets. Points \ensuremath{x} and \ensuremath{y} will be considered separated, by neighbourhoods, 
by closed neighbourhoods, by a continuous function, precisely by a function, 
iff their singleton sets $\{x\}$ and $\{y\}$ are separated according to the corresponding criterion.

Subsets A and \ensuremath{B} are {\em separated by neighbourhoods} iff
A and \ensuremath{B} have disjoint neighbourhoods, i.e. 
iff  
$ i_{AB} : \ensuremath{X} \longrightarrow  \{A \leftrightarrow  \ensuremath{x} \leftrightarrow  B\}$  factors as $$ 
X \longrightarrow   \{A \leftrightarrow  U_A {\small\searrow} \ensuremath{x} {\small\swarrow} U_B \leftrightarrow  B\} \longrightarrow  \{A \leftrightarrow  U_A = \ensuremath{x} = U_B \leftrightarrow  B\} $$
here the disjoint neighbourhoods of A and \ensuremath{B} are the preimages of open subsets ${A,U_A}$ and ${U_B,B}$ of 
$ \{A \leftrightarrow  U_A {\small\searrow} \ensuremath{x} {\small\swarrow} U_B \leftrightarrow  B\}$,  resp. 
They are {\em separated by closed neighbourhoods} 
iff they have disjoint closed neighbourhoods, i.e.
$i_{AB}$ factors as 
$$
X \longrightarrow   \{A \leftrightarrow  U_A {\small\searrow} U'_A {\small\swarrow} \ensuremath{x} {\small\searrow} U'_B {\small\swarrow} U_B \leftrightarrow  B\} \longrightarrow  \{A\leftrightarrow U_A=U'_A= \ensuremath{x} = U'_B=U_B\leftrightarrow B\} 
.$$
They are {\em separated by a continuous function} iff 
there exists a continuous function \ensuremath{f} from the space \ensuremath{X} to the real line $\RR$ such that $f(A)=0$  and $f(B)=1$,
i.e.
the map $i_{AB}$ factors as
$$
X \longrightarrow   \{0'\} \union [0,1] \union \{1'\}  \longrightarrow  \{A \leftrightarrow  \ensuremath{x} \leftrightarrow  B\}  
$$
where points $0',0$ and $1,1'$ are topologically indistinguishable, 
and $0'$ maps to $A$, and $1'$ maps to $B$, and $[0,1]$ maps to $x$. 
Finally, they are {\em precisely separated by a continuous function}
iff there exists a continuous function \ensuremath{f} from \ensuremath{X} to $\RR$ such that the preimage $f^{ - 1}(\{0\})= A$ and $f^{ - 1}(\{1\})=B$.
i.e.~iff $i_{AB}$ factors as $$ 
X \longrightarrow   [0,1]  \longrightarrow  \{A \leftrightarrow  \ensuremath{x} \leftrightarrow  B\} 
$$ where $0$ goes to point $A$ and $1$ goes to point $B$.

These conditions are given in order of increasing strength: 
Any two topologically distinguishable points must be distinct, and any two separated points must be topologically distinguishable. Any two separated sets must be disjoint, any two sets separated by neighbourhoods must be separated, and so on.

%For more on these conditions (including their use outside the separation axioms), see the articles Separated sets and Topological distinguishability.
%Main definitions

The definitions below all use essentially the preliminary definitions above.

In all of the following definitions, \ensuremath{X} is again a topological space.
\begin{itemize}
\item[]    \ensuremath{X} is T0, or Kolmogorov, if any two distinct points in \ensuremath{X} are topologically
distinguishable. (It will be a common theme among the separation axioms to have
one version of an axiom that requires T0 and one version that doesn't.)
As a formula, this is expressed as
   $$ \{x\leftrightarrow y\} \longrightarrow  \{x=y\} \,\rightthreetimes\,  \ensuremath{X} \longrightarrow  \{*\}$$

\item[]    \ensuremath{X} is R0, or symmetric, if any two topologically distinguishable points in X
are separated, i.e.
         $$\{x{\small\searrow}y\} \longrightarrow  \{x\leftrightarrow y\} \,\rightthreetimes\,  \ensuremath{X} \longrightarrow  \{*\}$$

\item[]    \ensuremath{X} is T1, or accessible or Frechet, if any two distinct points in \ensuremath{X} are
separated, i.e.
$$     \{x{\small\searrow}y\} \longrightarrow  \{x=y\} \,\rightthreetimes\,  \ensuremath{X} \longrightarrow  \{*\}  $$
 Thus, \ensuremath{X} is T1 if and only if it is both T0 and R0. (Although you may
say such things as "T1 space", "Frechet topology", and "Suppose that the
topological space \ensuremath{X} is Frechet", avoid saying "Frechet space" in this
context, since there is another entirely different notion of Frechet space in
functional analysis.)

\item[]    \ensuremath{X} is R1, or preregular, if any two topologically distinguishable points in
X are separated by neighbourhoods. Every R1 space is also R0.

\item[]    \ensuremath{X} is weak Hausdorff, if the image of every continuous map from a compact
Hausdorff space into \ensuremath{X} is closed. All weak Hausdorff spaces are T1, and all
Hausdorff spaces are weak Hausdorff.

\item[]    \ensuremath{X} is Hausdorff, or T2 or separated, if any two distinct points in \ensuremath{X} are
separated by neighbourhoods, i.e. 
$$    \{x,y\} \hookrightarrow  \ensuremath{X} \,\rightthreetimes\,  \{x{\small\searrow}X{\small\swarrow}y\} \longrightarrow  \{x=X=y\} $$
Thus, \ensuremath{X} is Hausdorff if and only if it is both T0
and R1. Every Hausdorff space is also T1.

\item[]    \ensuremath{X} is $T2\frac12$, or Urysohn, if any two distinct points in \ensuremath{X} are separated by
closed neighbourhoods, i.e. 
$$
 \{x,y\} \hookrightarrow  \ensuremath{X} \,\rightthreetimes\,  \{x{\small\searrow}x'{\small\swarrow}X{\small\searrow}y'{\small\swarrow}y\} \longrightarrow  \{x=x'=X=y'=y\}
$$
Every T2$\frac12$ space is also Hausdorff.

\item[]    \ensuremath{X} is completely Hausdorff, or completely T2, if any two distinct points in
X are separated by a continuous function, i.e. 
$$     \{x,y\} \hookrightarrow  \ensuremath{X} \,\rightthreetimes\,   [0,1]\longrightarrow \{*\}
$$     where  $\{x,y\} \hookrightarrow  X$  runs through all injective maps from the discrete two
point space $\{x,y\}$.

Every completely Hausdorff space is
also T2$\frac 12 ½$.

\item[]    \ensuremath{X} is regular if, given any point \ensuremath{x} and closed subset $F$ in \ensuremath{X} such that \ensuremath{x} does
not belong to $F$, they are separated by neighbourhoods, i.e. 
$$    \{x\} \longrightarrow  \ensuremath{X} \,\rightthreetimes\,  \{x{\small\searrow}X{\small\swarrow}U{\small\searrow}F\} \longrightarrow  \{x=X=U{\small\searrow}F\}
$$
(In fact, in a regular
space, any such \ensuremath{x} and\ensuremath{F} will also be separated by closed neighbourhoods.) Every
regular space is also R1.

\item[]    \ensuremath{X} is regular Hausdorff, or T3, if it is both T0 and regular.[1] Every
regular Hausdorff space is also $T2\frac12$.

\item[]    \ensuremath{X} is completely regular if, given any point \ensuremath{x} and closed set $F$ in \ensuremath{X} such
that \ensuremath{x} does not belong to $F$, they are separated by a continuous function, i.e. 
$$
      \{x\} \longrightarrow  \ensuremath{X} \,\rightthreetimes\,  [0,1] \cup \{F\} \longrightarrow  \{x{\small\searrow}F\}
$$    where points $F$ and $1$ are topologically indistinguishable, $[0,1]$ goes to $x$,
and $F$ goes to $F$.

Every
completely regular space is also regular.

\item[]    \ensuremath{X} is Tychonoff, or T3$\frac12$, completely T3, or completely regular Hausdorff, if
it is both T0 and completely regular.[2] Every Tychonoff space is both regular
Hausdorff and completely Hausdorff.

\item[]    \ensuremath{X} is normal if any two disjoint closed subsets of \ensuremath{X} are separated by
neighbourhoods, i.e.
$$   \emptyset \longrightarrow \ensuremath{X} \,\rightthreetimes\,  \{x{\small\swarrow}x'{\small\searrow}X{\small\swarrow}y'{\small\searrow}y\} \longrightarrow  \{x{\small\swarrow}x'=X=y'{\small\searrow}y\}
$$
 In fact, by Urysohn lemma a space is normal if and only if any two disjoint
closed sets can be separated by a continuous function, i.e.
 $$   \emptyset \longrightarrow  \ensuremath{X} \,\rightthreetimes\,  \{0'\} \cup [0,1] \cup \{1'\} \longrightarrow  \{0=0'{\small\searrow}x{\small\swarrow}1=1'\} $$ 
where points
$0',0$ and $1,1'$ are topologically indistinguishable,
            $[0,1]$ goes to $x$, and both $0,0'$ map to point $0=0'$,  and both $1,1'$ map to point
$1=1'$.

\item[]    \ensuremath{X} is normal Hausdorff, or T4, if it is both T1 and normal. Every normal
Hausdorff space is both Tychonoff and normal regular.

\item[]    \ensuremath{X} is completely normal if any two separated sets $A$ and $B$ are separated by
neighbourhoods $U\supset A$ and $V\supset B$
      such that $U$ and $V$ do not intersect, i.e.%????
 $$\emptyset \longrightarrow  \ensuremath{X} \,\rightthreetimes\,  \{X{\small\swarrow}A\leftrightarrow U{\small\searrow}U'{\small\swarrow}W{\small\searrow}V'{\small\swarrow}V\leftrightarrow B{\small\searrow}X\} \longrightarrow  \{U=U',V'=V\}$$
      Every completely normal space is also normal.
% $$
%      \emptyset \longrightarrow  \ensuremath{X} \,\rightthreetimes\,  \\ \{X{\small\swarrow}x\leftrightarrow U{\small\searrow}U'{\small\swarrow}W{\small\searrow}V'{\small\swarrow}V\leftrightarrow y{\small\searrow}X\} \longrightarrow 
%\{X{\small\swarrow}x\leftrightarrow U=U'{\small\swarrow}W{\small\searrow}V'=V\leftrightarrow y{\small\searrow}X\} $$
 %     or in short 

\item[]    \ensuremath{X} is perfectly normal if any two disjoint closed sets are precisely
separated by a continuous function, i.e. 
$$    \emptyset\longrightarrow \ensuremath{X} \,\rightthreetimes\,  [0,1]\longrightarrow \{0{\small\swarrow}X{\small\searrow}1\}
$$   where $(0,1)$ goes to the open point $X$, and $0$ goes to $0$, and $1$ goes to $1$.

Every perfectly normal space is also
completely normal.

\item[] \ensuremath{X} satisfies $T_d$ iff each point $x$ is the intersection of an open set $U$ and a closed set $Z$, $\{x\}=U\cap Z$, i.e.
each map the map $\{x\} \longrightarrow  X$ has the left lifting property wrt 
$$
\{ U\bar Z\bar x{\small\searrow}  {U} Zx{\small\searrow} \bar {U} Zx ,\ U\bar {Z} x\leftrightarrow  U\bar Z\bar x{\small\searrow} \bar U\bar {Z} x\leftrightarrow  \bar U\bar Z\bar x{\small\searrow} \bar {U} Zx\leftrightarrow  \bar {U} Z \bar x\}
$$
$$ \longrightarrow  
\{U\bar Z\bar {x} =  \bar U\bar Z\bar x= \bar {U} Z \bar {x} \leftrightarrow  
  {U} Zx=\bar {U} Zx=U\bar {Z} x= \bar U\bar {Z} x\}
$$
% $$
%\{x\} \longrightarrow  \ensuremath{X} \,\rightthreetimes\,    \{ U\bar \ensuremath{Z} x\leftrightarrow  U\bar Z\bar x{\small\searrow}  \ensuremath{U} Zx{\small\searrow} \bar \ensuremath{U} Zx\leftrightarrow  \bar \ensuremath{U} Z \bar \ensuremath{x} \  U\bar \ensuremath{Z} x\leftrightarrow  U\bar Z\bar x{\small\searrow} \bar U\bar \ensuremath{Z} x\leftrightarrow  \bar U\bar Z\bar x{\small\searrow} \bar \ensuremath{U} Zx\leftrightarrow  \bar \ensuremath{U} Z \bar x\} \longrightarrow   \{x\leftrightarrow \bar x\}
%$$

\item[]  \ensuremath{X} is extremally disconnected if the closure of every open subset of $X$ is open, i.e.
  $$\emptyset\longrightarrow \ensuremath{X} \,\rightthreetimes\,  \{U{\small\searrow}Z',Z{\small\swarrow}V\}\longrightarrow \{U{\small\searrow}Z'=Z{\small\swarrow}V\}$$
  or equivalently
   $$\emptyset\longrightarrow \ensuremath{X} \,\rightthreetimes\,  \{U{\small\searrow}Z',Z{\small\swarrow}V\}\longrightarrow \{Z'=Z\}$$ 

\end{itemize}

It is not clear if the property of being sober can be expressed as a lifting property. 

\newpage

%References:

\end{document}